\input amssym.def
\input amssym.tex

\let\bd=\bf 
\let\pd=\partial
\let\x=\times
\let\a=\alpha
\let\b=\beta
\let\g=\gamma
\let\Del=\triangledown
\def\cj#1{{{\overline{#1}}}}

\input plncs.cmm
\contribution{
Elementary Background in Elliptic Curves }
\author{Antonia W.~Bluher}
\address{National Security Agency, 9800 Savage Road,
Fort George G.\ Meade, MD  20755-6000}

\abstract
{
This paper gives additional background in algebraic geometry 
as an accompaniment to the article, ``Formal Groups, Elliptic Curves,
and some Theorems of Couveignes''.  Section 1 discusses the addition law
on elliptic curves, and Sections 2 and 3 explain about function fields,
uniformizers, and powerseries expansions with respect to a uniformizer.
}

\titlea{1}{   Elliptic curves}

If $K$ is a field then the 2-dimensional
projective space ${\bf P}^2(K)$ is the set of triples $(X,Y,Z)\in
K\times K\times K$ such that $X,Y,Z$ are not all zero, modulo the
equivalence $(X,Y,Z)\cong (\l X,\l Y,\l Z)$ for all $0\ne\lambda\in K$.
Three points $(X_i,Y_i,Z_i)$ ($i=1,2,3$)
in ${\bf P}^2(K)$ are said to be {\bd colinear} if
there exist $\lambda,\mu,\nu\in K$,
not all zero, such that $\l X_i+\mu Y_i+\nu Z_i=0$ for $i=1,2,3$.
In particular, if $Z_i\ne0$ for all $i$ then the colinearity of the
three points $(X_i,Y_i,Z_i)$ in ${\bf P}^2(K)$ is equivalent
to the colinearity of the points $(x_i,y_i)$ in the affine plane,
where $x_i=X_i/Z_i$ and $y_i=Y_i/Z_i$.  This follows from the equation
$\l x_i+\mu y_i+\nu=0$ for all $i$.

A {\bf Weierstrass equation} is a cubic equation $W(X,Y,Z)=0$, where
$$W(X,Y,Z)=Y^2Z+a_1XYZ+a_3YZ^2-(X^3+a_2X^2Z+a_4XZ^2+a_6Z^3),\eqno(1.1)$$
$a_i\in K$.  Since this equation is homogeneous, we can think of
its solutions as lying in ${\bf P}^2(K)$.
A {\bd singular point} is a point $P\in{\bf P}^2(\cj K)$ such that
$W(P)=\partial W/\partial X (P)=\partial W/\partial Y (P)
=\partial W/\partial Z (P)=0$, where $\cj K$ denotes
the algebraic closure of $K$.  It can be shown that the Weierstrass
equation has no singular points iff $\Delta\ne0$, where
$\Delta$ is a certain polynomial in the $a_i$; see [3], p.~46.
Let us assume $\Delta\ne0$.
Then the set of solutions to (1.1) in ${\bf P}^2(\cj K)$
is called an {\bd elliptic curve}.
The elements of the elliptic curve are called {\bd points}.
If $L$ is any field containing all the Weierstrass coefficients $a_i$
then
$E(L)$ is defined to be the set of solutions to (1.1) in
${\bf P}^2(L)$. The Weierstrass equation can often be simplified.
For example, in characteristic zero a linear change of coordinates
puts the Weierstrass equation into the form $Y^2Z=X^3+a_4XZ^2+a_6Z^3$.
See [3], Appendix~A.
 
The only point on the elliptic curve which intersects the line
$Z=0$ is the point $(0,1,0)$.  This is called the {\bd point at
infinity}
and is often denoted by $O_E$. In contrast, the ``finite'' points
of $E(K)$ are the points $(X,Y,1)$ which satisfy the above equation,
and it is customary to identify these points with the points $(X,Y)$ in
the affine plane $K\times K$ such that
$$Y^2+a_1XY+a_3Y=X^3+a_2X^2+a_4X+a_6.\eqno(1.2)$$
 
$E(K)$ can be made into an abelian group with additive identity equal to
$O_E$ by insisting that
any three colinear points on $E(K)$ (counting multiplicities)
sum to $O_E$.
The proof that this group law is associative can be found in
[2] using B\'ezout's Theorem or in [3] using the Picard group.
In this context, the phrase ``counting multiplicities''
means to count the multiplicities of the roots of the
cubic which results
when one variable  is eliminated from the Weierstrass
equation by using the equation of the intersecting line.
We give several examples which should make this clear.
\medskip
\titlec{ Example 1.1:  Lines through the origin.}
The equation of a line through $(0,1,0)$  has the form
$\lambda X+\nu Z=0$, where $\lambda,\nu\in K$ are not both zero.
If $\lambda=0$ then the line is $Z=0$, and the Weierstrass equation
becomes
simply $0=W(X,Y,0)=X^3$.  This has one solution $X=0$ with multiplicity
three.
Thus $(0,1,0)$ intersects $Z=0$ with multiplicity three.
This gives $3(0,1,0)=(0,1,0)$, which is consistent
with the fact that $(0,1,0)$ is the identity.
If $\lambda\ne0$ then $X=cZ$, where
$c=-\nu/\lambda$.  The Weierstrass equation becomes
$$W(cZ,Y,Z)=Y^2Z+(a_1c+a_3)YZ^2-(c^3+a_2c^2+a_4c+a_6)Z^3=0.$$
This factors as
$$Z(Y^2+(a_1c+a_3)YZ-(c^3+a_2c^2+a_4c+a_6)Z^2)=Z(Y-y_1Z)(Y-y_2Z)$$
where $y_1,y_2$ belong to a quadratic extension of $K$ and satisfy
$y_1+y_2=-(a_1c+a_3)$, $y_1y_2=-(c^3+a_2c^2+a_4c+a_6)$.  The
intersection
of $E$ with the line is $\{(0,1,0),(c,y_1,1),(c,y_2,1)\}$, hence
$(0,1,0)+(c,y_1,1)+(c,y_2,1)=(0,1,0)$. In other words, if $(c,y_1,1)$
lies
on $E$ then
$-(c,y_1,1)=(c,y_2,1)$, where $y_1+y_2=-(a_1c+a_3)$.  This gives the
inversion formula:
$$-(x,y,1)=(x,-(y+a_1x+a_3),1).\eqno(1.3)$$
If it happens that $y=-(y+a_1x+a_3)$, that is, $2y=-(a_1x+a_3)$, then
$(x,y,1)$ is a two-torsion point.

\titlec{Example 1.2: Lines which miss the origin.}
The equation of a line which misses $(0,1,0)$ is
$\lambda X+\mu Y + \nu Z=0$, where $\mu\ne0$.  Since
$(0,1,0)$ is the only point on $E$ which intersects $Z=0$, the
intersection points
of the line with the elliptic curve all have a nonzero $Z$-coordinate.
Let $m=-\lambda/\mu$, $b=-\nu/\mu$, $y=Y/Z$, $x=X/Z$. Then the three
intersection points $(x_i,y_i,1)$, $i=1,2,3$,
satisfy the Weierstrass equation and $y_i=mx_i+b$.
Another way to view this is: suppose $(x_1,y_1,1)$ and $(x_2,y_2,1)$ are
two
points on the elliptic curve such that $x_1\ne x_2$, and define
$m=(y_1-y_2)/(x_1-x_2)$, $b=y_1-mx_1$.
Substitute $y=mx+b$ into
the affine form of the Weierstrass equation to get a cubic equation
in $x$
of the form $x^3+Bx^2+Cx+D=0$, where the coefficients explicitly depend
on the $a_i$ and on $m$ and $b$.  Explicitly, $B=-m^2-a_1m+a_2$.
Two roots of this cubic equation
are $x_1$ and $x_2$.  Let $x_3$ denote the third root. Then
$$x^3+Bx^2+Cx+D=(x-x_1)(x-x_2)(x-x_3)=x^3-(x_1+x_2+x_3)x^2+(\cdots)x+(\cdots),$$so that $x_3=-x_1-x_2-B.$  This yields the addition formula
$$(x_1,y_1,1)+(x_2,y_2,1)=(x_3,-(y_3+a_1x_3+a_3),1)
\qquad{\rm if\ }x_1\ne x_2,$$
$$x_3=-x_1-x_2+m^2+a_1m-a_2,\qquad m=(y_1-y_2)/(x_1-x_2),$$
$$y_3=m x_3 + y_1-mx_1.$$
\titlec{ Example 1.3: Duplication formula.}
How do we compute $P+P$?  Assume $P=(x_1,y_1,1)$ is not
a two-torsion point, that is, $2y_1\ne-(a_1x+a_3)$.  Then $-(P+P)$
is the point
$(x_2,y_2,1)$ on the curve such that $2(x_1,y_1,1)+(x_2,y_2,1)=(0,1,0)$.
Let $\l X +\mu Y +\nu Z$ be the line which
passes through these points with correct multiplicity.  It misses
the origin, so $\mu\ne 0$.  Thus the line has the form
$y-y_1=m(x-x_1)$, where $y=Y/Z$, $x=X/Z$.
First we will find $m$,
then we will find $x_2,y_2$.
When $y$ is replaced by $m(x-x_1)+y_1$ in the Weierstrass equation,
we get a cubic  $x^3+Bx^2+Cx+D$,
and it should equal $(x-x_1)^2(x-x_2)$.
Here $B=a_2-m^2-a_1m$, $C=m^2(2x_1)+m(-2y_1+a_1x_1-a_3)+(a_4-a_1y_1)$.
If we differentiate the cubic and set $x=x_1$ we are
supposed to get zero.
Hence $3x_1^2+2Bx_1+C=0$.  We can solve for $m$ in this relation;
in fact $m$ will just be the slope of the tangent line to the curve
$E$ at $P$, namely $m=(dy/dx)(P)$.  It turns out
$$m=(3x_1^2+2a_2x_1+a_4-a_1y_1)/(a_1x_1+2y_1+a_3).$$
Since $(x-x_1)^2(x-x_2)=x^3+Bx^2+Cx+D$,
$$x_2=-2x_1-B=-2x_1-a_2+m^2+a_1m,\qquad y_2=m(x_2-x_1)+y_1.$$
Finally,
$$2(x_1,y_1,1)=-(x_2,y_2,1)=(x_2,-(y_2+a_1x_2+a_3),1).$$

\titlea{2}  { Function fields, local rings, and uniformizers}

Let $K$ be any field, and
let $C=C(X,Y,Z)$ be a curve in ${\bf P}^2(\cj K)$.
This means that $C$ is an irreducible homogeneous polynomial
with coordinates in $\cj K$.  Assume the coefficients of $C$
belong to $K$, and let $C(K)=\{\,P\in{\bf P}^2(K)\,|\,C(P)=0\,\}$.
A point in $C(\cj K)$ is called
{\bd nonsingular} if $\triangledown_P C\ne0$, where
$$\Del_PC=\left({\pd\over\pd X}C(P),{\pd\over\pd Y}C(P),
{\pd\over\pd Z}C(P)\right).$$
The {\bd function field} $K(C)$ of $C$ over $K$ is the set of all
quotients
$$p(X,Y,Z)/q(X,Y,Z)$$, where $p$ and $q$ are homogeneous polynomials of
the same degree with coefficients in $K$ such that $C$ does not
divide $q$, modulo the equivalence: $p/q\cong p'/q'$ iff
$pq'-qp'$ is divisible by $C$.  The addition and multiplication
in $K(C)$ are defined in the same way as for rational functions.
Equivalently, if we set $x=X/Z$ and $y=Y/Z$ then $K(C)$
is the quotient field of the integral domain
$$K[x,y]/(C(x,y,1)).$$
This field has transcendence degree~1, thus any two elements of the
function field satisfy some algebraic relation.
Notice that $p(\l X,\l Y,\l Z)/q(\l X,\l Y,\l Z)
=p(X,Y,Z)/q(X,Y,Z)$ because
of the homogeneity, so $p(P)/q(P)$ makes sense for $P\in C(\cj K)$
provided $q(P)\ne0$.
A function $f\in K(C)$ is said to be {\bd defined} at a point
$P\in C(\cj K)$ if
there exists $p/q$ in the equivalence class of $f$ such that
$q(P)\ne 0$.
If $f$ is defined at $P$ then the {\bd value} of $f$ at $P$,
denoted $f(P)$,  is defined as
$p(P)/q(P)$ for any (hence every) $p/q$ which is in the
equivalence class
of $f$ and for which $q(P)\ne0$.  If $P\in C(K)$ and $f\in K(C)$
is defined at $P$ then $f(P)\in K$.
 
As an example, let us show that for the Weierstrass equation,
$Z/X$ is defined at $(0,1,0)$.
Since $Z(Y^2+a_3 YZ-a_6Z^2)=X(-a_1 YZ+X^2+a_2XZ+a_4Z^2)$,
$$Z/X={-a_1 YZ+X^2+a_2XZ+a_4Z^2\over Y^2+a_3 YZ-a_6Z^2}.$$
Now $Y^2+a_3 YZ-a_6Z^2$ does not vanish at $(0,1,0)$, so $Z/X$ is indeed
defined at $(0,1,0)$.
 
{\it For the remainder of this section assume $P\in C(K)$.}
Let $\Omega_P$ be the ring of functions in $K(C)$ which are
defined at $P$, and let $M_P$ be the ideal of $\Omega_P$ consisting
of functions which take the value zero at $P$.  $\Omega_P$ is called
the {\bf local ring of~$C$ at~$P$}. If we identify
the set of constant functions with $K$ then
$$\Omega_P=K\oplus M_P\eqno(3.1)$$
(internal direct sum) because $f=f(P)+(f-f(P))$.
 
\lemma {2.1} 
{
$\Omega_P-M_P=\Omega_P^\x$. $M_P$ is the
unique maximal ideal of $\Omega_P$.
}
 
\proof{}
If $f\in\Omega_P-M_P$
then $f$ can be written as $F(X,Y,Z)/G(X,Y,Z)$,
where $F,G$ are homogeneous polynomials of the same degree and $F,G$
do not vanish at $P$. Since $1/f$ can be written $G/F$,
$1/f\in\Omega_P$. This shows that $\Omega_P-M_P\subset\Omega_P^\x$.
The reverse inclusion is also true, for if $f\in\Omega_P^\x$
then there is $g\in\Omega_P$ with $fg=1$. Then $f(P)g(P)=1$,
so $f(P)$ is a unit in $K$.  In particular, $f$ is not in $M_P$.
We have proved the first statement. For the second statement,
let $I$ be any proper ideal of $\Omega_P$. Then $I\cap\Omega_P^\x=
\emptyset$; otherwise $I$ would equal $\Omega_P$.  Thus
$I$ misses the complement of $M_P$ completely; equivalently,
$I\subset M_P$.\qed
 
Our next goal is to show that if $P$ is a nonsingular point
(that is, $\Del_PC\ne0$) then $M_P$ is principal. First we need
some lemmas. These lemmas are true in more generality
than we state them.
 
\lemma {2.2 (Nakayama's Lemma)}
{
If $A$ is a finitely
generated $\Omega_P$-module and $M_PA=A$ then $A=\{\,0\,\}$.
}
 
\proof{}  Let $u_1,\ldots,u_n$ be a set of generators for $A$
with $n$ as small as possible. Since $u_1\in A=M_PA$, there exist
$\mu_i\in M_P$ such that $u_1=\sum \mu_i u_i$.  Now $1-\mu_1\in
\Omega_P^\x$, so we can solve for $u_1$ in terms of $u_2,\ldots,
u_n$. This contradicts the minimality of $n$.\qed
 
\corollary {2.3}
{
$\cap_{n=1}^\infty M_P^n=\{\,0\,\}$.  Thus
$M_P\supsetneq M_P^2$.
}
 
\proof{} Let $A=\cap M_P^n$. $A$ is finitely generated by the Hilbert
Basis Theorem [1], p.~13, so $A=\{\,0\,\}$ by Nakayama's Lemma.\qed
 
\lemma {2.4}
{
If $P$ is a nonsingular point of $C$ then
$M_P/M_P^2$ is isomorphic to $K$.
}
 
\proof{} Observe that $M_P/M_P^2$ is an $\Omega_P/M_P$-vector space,
and $\Omega_P/M_P$ is canonically isomorphic to $K$ by (3.1).
We just have to prove this vector space is one-dimensional,
or equivalently that its dual is one-dimensional.
So let $\l:M_P/M_P^2\to K$ be a $K$-linear map.
Since $\Omega_P=K+M_P$, we can think of $\l$ as a $K$-linear map
from $\Omega_P$ into $K$ which is trivial on $M_P^2+K$.
 
We will construct a $K$-linear map $\vartheta:\Omega_P\to
K$ which is trivial on $K+M_P^2$, and we will prove that
$\lambda$ is proportional to $\vartheta$.
Fix coordinates for $P$: $P=(X_0,Y_0,Z_0)$. Fix $R\in K^3$
such that $R\cdot P=1$. Let $C(X,Y,Z)=0$ be the equation
of the curve. Let
$\Del_P C=(\pd C/\pd X(P),\pd C/\pd Y(P),\pd C/\pd Z(P))$,
the gradient of $C$ at $P$. We claim $\Del_P C$ and $R$
are linearly independent. First observe $\Del_P C\ne0$
since the curve is nonsingular at $P$ by hypothesis.
By Euler's identity $(X\pd/\pd X+Y\pd/\pd Y + Z \pd/\pd Z)C=
{\rm deg}(C) C$.  (The proof of Euler's identity is that
for any monomial $X^aY^bZ^c$, we have
$(X\pd/\pd X+Y\pd/\pd Y + Z \pd/\pd Z)(X^aY^bZ^c)=(a+b+c)X^aY^bZ^c$.)
Evaluating Euler's identity at $P$ gives
$P\cdot\Del_P C=deg(C)\,C(P)=0$.
Since $P\cdot R\ne0$, this shows $\Del_P C$ is linearly independent
from $R$.  It follows that the space of vectors which are orthogonal
to $R$ and $\Del_PC$ is one-dimensional, spanned by a
vector $T$.
 
Define $\vartheta:\Omega_P\to K$ by
$$\vartheta(f)=T\cdot\Del_P(f).$$
We claim $\vartheta$ is well-defined, $K$-linear, and  surjective.
To prove it is well-defined, it suffices to show
$T\cdot\Del_P(F/G)=0$ if $F$ vanishes identically on the curve
and $G(P)\ne0$. In that case, $C$ divides $F$. Let $h=F/(CG)$.
Then
$$T\cdot\Del_P (F/G)=C(P) T\cdot\Del_P h+h(P)T\cdot\Del_P C
=0+0$$
as required. The map $\vartheta$
is certainly $K$-linear.
If $S$ is any vector in $K^3$ such that $S\cdot P=0$ but
$S\cdot T\ne0$ then
$$T\cdot\Del_P (S\cdot(X,Y,Z)/R\cdot(X,Y,Z))
={T\cdot S\over R\cdot P}
-{(T\cdot R)(S\cdot P)\over (R\cdot P)^2}
={T\cdot S\over R\cdot P} \ne0,\eqno(3.1)$$
so the map is nontrivial. Any nontrivial
$K$-linear map into $K$
is surjective.  Let us show that the kernel contains
$K+M_P^2$.  Certainly $T\cdot\Del_P$ annihilates the constants
$K$.  If $f,g\in M_P$ then $\Del_P(fg)=0$ by the product
rule, so $M_P^2$ is in the kernel also.

It remains to show that an arbitrary linear map $\l:\Omega_P\to K$
which is trivial on $K+M_P^2$ coincides with a multiple of $\vartheta$.
The map $\lambda$ obeys a product rule, for if $f=f(P)+f_1$,
$g=g(P)+g_1$ where $f,g\in\Omega_P$ then $f_1,g_1\in M_P$,
and $fg=f(P)g+g(P)f-f(P)g(P)+f_1g_1$, so
$$\l(fg)=f(P)\l(g)+g(P)\l(f).$$
It follows that $\l$ is completely determined by its
values at $X/\rho$, $Y/\rho$, and $Z/\rho$,
where $\rho=R\cdot(X,Y,Z)$.
Denote these by values by $\a,\b,\g$.  We claim $(\a,\b,\g)$ is
proportional to $T$. This will imply there is a one-dimensional
space of such $\lambda$; in particular
$\l$ must be proportional to $\vartheta$.
We must show $(\a,\b,\g)$ is orthogonal to both $R$ and
$\Del_PC$.  Let $R=(R_1,R_2,R_3)$ and $\rho=R_1X+R_2Y+R_3Z$. Then
$$\eqalign{R\cdot(\a,\b,\g)&=R_1\l(X/\rho)
+R_2\l(Y/\rho)+R_3\l(Z/\rho)\cr
&=\l((R_1X+R_2Y+R_3Z)/\rho)=\l(1)=0,\cr}$$
$$\Del_P C\cdot(\a,\b,\g)=\l(C(X/\rho,Y/\rho,Z/\rho))=\l(C(X,Y,Z)/
\rho^{\deg (C)})
=0.$$
\qed
 
\corollary {2.5}
{
Supose $P$ is nonsingular.
Let $u\in M_P$, and let $R,T$ be
nonzero vectors in $K^3$ such that $R\cdot P\ne0$, $T\cdot R=0$,
and $T\cdot \Del_PC=0$. Then $u\in M_P-M_P^2$ iff
$T\cdot\Del_Pu\ne0$. Also, if $S\cdot P=0$ but $S\cdot T\ne 0$ then
$S\cdot(X,Y,Z)/R\cdot(X,Y,Z)$ belongs to $M_P-M_P^2$.
}
 
\proof{} In the notation of the preceding proof, the hypothesis
states that $\vartheta(u)\ne0$. Since $\vartheta$ induces an
isomorphism of $M_P/M_P^2$ with $K$, we know $u\not\in M_P^2$
iff $\vartheta(u)\ne0$.  The last statement follows from (3.1).
\qed
 
\proposition {2.6} 
{
Suppose $P$ is nonsingular.
Then $M_P$ is a principal ideal of $\Omega_P$.
In fact, if $u\in M_P$ then $M_P=u\Omega_P\iff u\not\in M_P^2$.
}
 
\proof{} $M_P\ne M_P^2$ by Lemma 2.2. Let $u\in M_P-M_P^2$.
Then
$$K\cong{\Omega_P\over M_P}\cong{u\Omega_P\over uM_P}\to
{u\Omega_P\over M_P^2}.$$
The image of 1 in this composite map is $u\bmod M_P^2$,
which is nonzero. Thus the kernel of the
projection map from $u\Omega_P/uM_P$ into $u\Omega_P/M_P^2$
has to be trivial. This shows $uM_P=M_P^2$.
 
Also, $M_P=uK+M_P^2$ because $M_P/M_P^2$ is a one-dimensional
$K$-vector space by Lemma~3.4. Thus
$$u\Omega_P=u(K+M_P)=uK+uM_P=uK+M_P^2=M_P.$$
        \qed
 
\corollary {2.7}
{  
Suppose $P$ is nonsingular.
Let $0\ne f\in K(C)$, $u\in M_P-M_P^2$.
There is a unique integer $v(f)$ such that $f=u^{v(f)}g$
with $g\in\Omega_P^\x$. This integer does not depend on the choice
of $u$.
}
 
\proof{}  Write $f=F/G$ with $F,G$ homogeneous of degree $m$.
Choose $R\in K^3$ so that $P\cdot R\ne0$. Let
$\rho$ be the homogeneous polynomial $\rho=R\cdot(X,Y,Z)$.
Then $f=f_1/f_2$, where $f_1=F/\rho^m$ and $f_2=G/\rho^m$.
Note that $f_1,f_2\in\Omega_P$.  By Corollary~3.3,
$$\Omega_P-\{0\}=\sqcup_{n=0}^\infty (M_P^n-M_P^{n+1}),$$
where $M_P^0=\Omega_P$.  Clearly $M_P^{n}-M_P^{n+1}=
u^n \Omega_P^\x$, since $M=u\Omega_P$. Define $n_1,n_2$ by
$f_i\in M_P^{n_i}-M_P^{n_i+1}$, and let $n=n_1-n_2$.
Then $fu^{-n}=f_1u^{-n_1}/(f_2 u^{-n_2})\in\Omega_P^\x$.
\qed
 
$\Omega_P$ is called the {\bf local ring} at $P$. If
$P$ is nonsingular, then a {\bf uniformizer}
is an element of $M_P$ such that $M_P=u\Omega_P$;
equivalently it is an element of $M_P-M_P^2$.
The integer $v(f)$ is called the {\bf valuation} of $f$
at $P$.   If $v(f)>0$, $f$ has a {\bf zero of order $v(f)$ at
$P$}; if $v(f)<0$, $f$ has a {\bf pole of order $-v(f)$ at $P$}.
If $v(f)=0$, we say $f$ is {\bf finite and nonvanishing at $P$}.

\titlec{ Example 2.8} Let $C=W$ be given by the
Weierstrass equation (1.1).
Let us show that $X/Y$ is a uniformizer at
$(0,1,0)$.  In Corollary 2.5 we take $R=P=(0,1,0)$ and
compute from (1.1) that $\Del_PW=(0,0,1)$. Thus $T=(1,0,0)$.
Now $X/Y$ vanishes at $P$, and $T\cdot\Del_P(X/Y)=\pd/\pd X(X/Y)|_P=1
\ne0$. So $X/Y$ is a uniformizer by Corollary~2.5.
\qed
\titlec{ Example 2.9} Let us compute $v(Z/X)$ at the point
$P=(0,1,0)$ on an elliptic curve. By (1.1), $ZV=X^3$, with
$$V=Y^2+a_1XY+a_3YZ-a_2X^2-a_4XZ-a_6Z^2.$$
Then $Z/X=X^2/V=(X/Y)^2(Y^2/V)$.
Since $Y^2/V\in\Omega_P^\x$ and $X/Y$ is a uniformizer,
we see that $v(Z/X)=2$.  Thus $x=X/Z$ has a pole of order~2
at the origin. Since $Y/Z=(Y/X)(X/Z)$, $y=Y/Z$ has a pole
of order~3 at the origin.
\qed
 
\titlea{3} {Power series expansions}

\lemma {3.1}
{
Let $C=C(X,Y,Z)$ be a projective curve
over $\cj K$ which is nonsingular
at $P$. Suppose $C$ has coefficients in $K$ and $P\in C(K)$.
Let $u\in K(C)$ be a uniformizer at $P$.
If $f$ is a nonzero function in $K(C)$ and
$v(f)=a$ then there exist unique constants
$f_j\in K$ for $j\ge a$ such that for any $N\ge a$,
$$f
-\sum_{j=a}^N f_ju^j\qquad{\rm has\ a\ zero\ of\ order\ at\ least\ }
N+1.$$
Moreover, $f_a\ne0$.
}
 
\proof{}  Let $g=u^{-a}f\in\Omega_P^\x$, $f_a=g(P)$.
This is the only constant for which $g-f_a\in M_P$,
so it is the only constant for which $v(f-u^af_a)=v(u^a(g-f_a))>a$.
>From this observation, the lemma may be easily proved by induction.
\qed
 
\corollary {3.2}
{
A choice of uniformizer $u$ at $P$ determines
a one to one homomorphism of rings
$$\Psi_u:\Omega_P\hookrightarrow K[[\tau]].$$
If $0\ne V=\sum v_i\tau^i\in K[[\tau]]$, define $\deg(V)$
to be the smallest integer $n$ such that $v_n\ne0$, and
define $\deg(0)=\infty$.
If one puts metrics on $\Omega_P$ and on $ K[[\tau]]$ by
the rules $|f|=c^{v_P(f)}$ for $f\in\Omega_P$,
$|V|=c^{\deg(V)}$ for $V\in  K[[\tau]]$, where $0<c<1$,
then $\Psi_u$ is an isometry.  In particular, $\Psi_u$ is
continuous with respect to the topologies on $\Omega_P$,
$ K[[\tau]]$ which are induced by the above metrics.
The image of $\Psi_u$ is dense in $ K[[\tau]]$.
}
 
\proof{}  The homomorphism $\Psi_u$ is defined by:
$\Psi_u(f)=\sum f_j\tau^j$, where the $f_j$ are as in Lemma~3.1.
It is easy to see (using the uniqueness of the coefficients $f_j$)
that $\Psi_u$ is a ring homomorphism.
$\Psi_u$ is an isometry because $v_P(f)=a$ implies
$f_a$ is the first nonzero coefficient. If $\Psi_u(f)=0$ then
$|f|=|\Psi_u(f)|=|0|=0$, so $f\in\cap M_P^n$.
By Corollary~2.3, $f$ must be zero; thus $\Psi_u$ is one to one.
If $g$ is any polynomial with coefficients in $ K$
then $g(u)\in\Omega_P$ and
$\Psi_u(g(u))=g(\tau)$. This shows $\Psi_u(\Omega_P)$ is dense in
$ K[[\tau]]$.\qed

\begref{References}{3.}

\refno {1.} W.~Fulton, {\it Algebraic Curves}, Addison-Wesley, 1969
\refno {2.} A.~W.~Knapp, {\it Elliptic Curves}, {\sl Math.~Notes\ }{\bf40},
Princeton University Press, Princeton, New Jersey, 1992
\refno {2.} J. H. Silverman, {\it The Arithmetic of Elliptic Curves},
Springer-Verlag, New York, 1986

\end